\documentclass[12pt]{article}
\usepackage{subeqn,multirow,url}
\usepackage{color}

\title{\bf On Pythagorean triplets}

\author{\bf Palash B. Pal \\
Department of Physics, University of Calcutta,\\ 92 APC Road, Calcutta
700009, India}

\date{}

\evensidemargin=\oddsidemargin
\newtheorem{notn}{Notation}
\newtheorem{theorem}{Theorem}
\def\proof{\par \medskip\noindent{\bf Proof:\ }\rm\hangindent=5mm
  \parindent=10mm}

\makeatletter
\let\c@notn\c@theorem
\makeatother

\def\Eqn#1{Eq.\,(\ref{#1})}
\def\theo#1{Theorem \ref{#1}}
\def\tabl#1{Table\,\ref{#1}}
\def\mod{\mathrel{\rm mod}}
\def\bsup#1{^{(#1)}}

\def\trip#1 #2 #3 {\left(\displaystyle{#1 \;\; #2 \atop
    #3}\right)\vphantom{{a^A \over a^A} \over
    {b^2 \over b^2}}}

\def\trip#1 #2 #3 {\left(\displaystyle{#3 \atop
    #1 \quad #2}\right)\vphantom{{a^A \over a^A} \over
    {b^2 \over b^2}}
}

\def\options(#1){\left\{ \begin{array}{c} #1 
  \end{array} \right\}}

\begin{document}

\maketitle

\begin{abstract}

  We discuss properties of diophantine solutions of the Pythagoras
  equation, $a^2+b^2=c^2$, where the three numbers have no common
  factor.  Some of the highlights are: (1) All triplets for which $c$
  (called the `peak') is non-prime can be deduced from the triplets
  with prime peaks; (2) If a peak has $n+1$ prime factors, there are
  $2^n$ independent solutions of the Pythagoras equation; (3) All
  Pythagorean peaks have to be of the form $12k+1$ or $12k+5$ for
  integer $k$; (4) A Pythagorean peak cannot have 3, or any number of
  the form $12k+7$ or $12k+11$, as its prime factors.

\end{abstract}

\centerline{------------------------------------------------------------}

\bigskip

Pythagorean triplets are defined to be integer valued solutions of the
equation
\begin{eqnarray}
  a^2 + b^2 = c^2 .
  \label{p}
\end{eqnarray}
There are infinite number of solutions.  The question that we want to
ask is this: can all such solutions be generated from a smaller number
of basic solutions?  By `generate', we mean that all three numbers in
the triplet can be obtained by simple arithmetical operations from
some other triplet, without having to take recourse of any search
program of any sort.  If a triplet cannot be generated from other
triplets, we will call it a {\em basic triplet}.

We can of course disregard solutions with negative integers when we
look for basic triplets.  Any solution with positive integers will be
accompanied by other solutions in which one or more of the three
numbers appearing in \Eqn{p} are replaced by their negatives.  Thus,
in order to find the basic solutions, we can take
\begin{eqnarray}
  a>0, \qquad  b>0, \qquad c>0.
  \label{+ve}
\end{eqnarray}
In this case, clearly
\begin{eqnarray}
  a,b < c \,.
\end{eqnarray}
Also,
\begin{eqnarray}
  a+b > c \,.
  \label{treq}
\end{eqnarray}
Although \Eqn{treq} can be seen as the triangle inequality, we
emphasize that we need not think about any triangle for our entire
discussion.  We can just treat \Eqn{p} as a diophantine equation, and
refer to $c$ not as the hypotenuse but as the {\em Pythagorean peak}.
The other two numbers will be referred to as the {\em Pythagorean
  base}.  Even in this approach, \Eqn{treq} follows from the fact that
$(a+b)^2 = c^2+2ab > c^2$.  In what follows, all algebraic symbols
would refer to non-negative integers, a restriction that will not be
mentioned every time.

Next, from every solution we can construct other solutions by
multiplying each of the numbers by the same integral factor.  We do
not need to consider those other solutions.  Thus, we need to look
only at solutions for which
\begin{eqnarray}
  \gcd (a,b,c) = 1 \,.
  \label{gcd1}
\end{eqnarray}
Henceforth, we will only discuss solutions of \Eqn{p} satisfying
\Eqn{+ve} and \Eqn{gcd1}.  Only these solutions will be called
{\em Pythagorean triplets}.

\begin{theorem}\label{th:2odd}

  Two of the numbers $a,b,c$ are odd, the other is even.

  \proof The number of odd numbers must be even in order to satisfy
  \Eqn{p}.  The number cannot be zero, because then all three numbers
  will be even, violating \Eqn{gcd1}.  Hence, it must be 2.

\end{theorem}

\begin{theorem}\label{oddc}

  The number $c$ must be odd.

  \proof Consider the contrary of the statement, i.e., $c$ is even.
  Then both $a$ and $b$ are odd, according to \theo{th:2odd}.  So let
  us say $a=2n_a+1$, $b=2n_b+1$ with two integers $n_a$ and $n_b$.
  Then
  \begin{eqnarray}
    a^2+b^2 = (2n_a+1)^2 + (2n_b+1)^2 = 2 \mod 4.
  \end{eqnarray}
  But since $c$ is even, we have $c^2 = 0 \mod 4$.  This is a
  contradiction.  This proves that $c$ must be odd.

\end{theorem}

One of the two numbers $a$ or $b$ must therefore be even, and the
other odd.

\begin{notn}
  Henceforth, we will name them in such a way that $a$ is even and $b$
  is odd.
\end{notn}

There are certainly solutions to \Eqn{p} of the form
\begin{eqnarray}
  a = 2pq, \qquad b = |p^2-q^2|, \qquad c = p^2+q^2 
  \label{pq}
\end{eqnarray}
for integers $p$ and $q$.  Some features of such solutions are worth
discussing.

\begin{theorem}

  Among the two numbers $p$ and $q$, one must be even and one odd.

  \proof If both are even or both are odd, $c$ cannot be odd, as
  required from \theo{oddc}.

\end{theorem}

\begin{notn}
  We will take $p$ to be even and $q$ to be odd.  This is no loss of
  generality. 

\end{notn}

Hence we are led to the conclusion that, for solutions of the type
given in \Eqn{pq}, the number $a$ must be divisible by 4.

Is the converse true?  Certainly, if $a$ is divisible by 4, we can
write it as $2pq$, with an even $p$ and odd $q$, or vice versa.  And
then, if we construct $b$ and $c$ as shown in \Eqn{pq}, they will
satisfy Pythagoras' equation.

So then the question is: are all solutions to \Eqn{p} of the form
given in \Eqn{pq}?  The answer is yes.

\begin{theorem}\label{th:allsol}

  All solutions of \Eqn{p} can be expressed in the form of
  \Eqn{pq} with the help of an even and an odd integer.

  \proof As already said, a solution of the form of \Eqn{pq} is
  equivalent to having $a=4n_a$ for some integer $n_a$.  So the
  question boils down to demonstrating that there cannot be any
  solution of \Eqn{p} in which $a$ is not divisible by 4.

  Since $a$ is even, the only other possibility to be explored is that
  of $a=4n_a+2$.  If $a$ is of this form, then
  \begin{eqnarray}
    a^2 = 16n_a^2 + 16n_a + 4 = 4 \mod 8.
    \label{asqmod8}
  \end{eqnarray}
  However,
  \begin{eqnarray}
    c^2-b^2 = (2n_c+1)^2 - (2n_b+1)^2 = 4 \Big( n_c(n_c+1) -
    n_b(n_b+1) \Big) .
    \label{bcmod8}
  \end{eqnarray}
  Any product of the form $n(n+1)$ is even, so the number within the
  parentheses is even.  Hence the right side of \Eqn{bcmod8} must be
  divisible by 8, in contradiction with \Eqn{asqmod8}.  This proves the
  theorem.   

\end{theorem}

\begin{table}
  \caption{List of Pythagorean triplets generated from \Eqn{pq} with
    $p,q<10$.  The entries marked with a cross indicate that the
    corresponding $p$ and $q$ are not relatively prime.}\label{list}

  $$
  \begin{array}{c|cccccccc}
      q \rightarrow & 1 & 3 & 5 & 7 & 9 \\
      p \downarrow \\ 
      \hline 
      2 &  \trip 4 3 5 & \trip 12 5 13 &
      \trip 20 21 29 & \trip 28 45 53 & \trip
      36 77 85 \\
      4 & \trip 8 15 17 & \trip 24 7 25 & \trip 40 9 41 & \trip 56 33
      65 & \trip 72 65 97 \\
      6 & \trip 12 35 37 & \times & \trip 60 11 61 & \trip 84 13
      85 & \times \\
      8 & \trip 16 63 65 & \trip 48 55 73 & \trip 80 39 89 & \trip 112
      15 113 & \trip 144 17 145 \\
    \end{array}
  $$
    
\end{table}
Clearly, if $p$ and $q$ contain a common factor $r$, then all numbers
generated by \Eqn{pq} would contain a factor $r^2$.  Therefore, the
recipe for producing Pythgorean triplets is simple: take an even
integer $p$ and an odd integer $q$ which are relatively prime, and
construct the numbers $a$, $b$, $c$ as given in \Eqn{pq}.  This recipe
will produce all Pythagorean triplets.

Fine.  But we ask ourselves the question: instead of taking recourse
of $p$ and $q$ to obtain the Pythagorean triplets, can we generate
some Pythagorean triplets from the others, and bootstrap all of them
starting from a few?

To motivate the road to an answer, let us rewrite the Pythagoras
equation in the form
\begin{eqnarray}
  \Bigg( {a \over c} \Bigg)^2 + \Bigg( {b \over c} \Bigg)^2 = 1.
\end{eqnarray}
This shows that $a/c$ and $b/c$ can be coordinates of a unit circle
\cite{stewart}.  The ratios are rational numbers.  Thus, Pythagorean
triplets can be identified by all points on the unit circle for which
both $x$ and $y$ coodinates are rational.  Let us call these points on
the unit circle as {\em Pythagorean points}.

We now define an angle $\theta$ by the relation
\begin{eqnarray}
  \tan\theta = {a \over b}.
\end{eqnarray}
For any Pythagorean point, $\tan\theta$ will be rational.  Certainly,
all ratios involving $\tan\theta$ and its powers will also be
rational.  In particular, consider
\begin{eqnarray}
  \tan2\theta = {2\tan\theta \over 1 - \tan^2\theta} = {2ab \over
    b^2-a^2} \,. 
  \label{tan2theta}
\end{eqnarray}
This says that $2ab$ and $|b^2-a^2|$ can constitute a Pythagorean
base.  This is no big news, since these form a subset of the numbers
defined in \Eqn{pq}.

But we can do more.  Consider
\begin{eqnarray}
  \tan3\theta = {3\tan\theta - \tan^3\theta \over 1 - 3\tan^2\theta} =
              {3ab^2 - a^3 \over b^3 - 3a^2b} \,.
\end{eqnarray}
This shows that we can consider $3ab^2 - a^3$ and $b^3 - 3a^2b$ (or,
their aboslute values if any of them is negative) as a Pythagorean
base as well.  Indeed, it is easy to see that
\begin{eqnarray}
  (3ab^2 - a^3)^2 + (b^3 - 3a^2b)^2 = (a^2+b^2)^3 \,.
\end{eqnarray}
The right side is not the square of an integer in general, but if $a$
and $b$ satisfy \Eqn{p}, the right side can be written as $(c^3)^2$.
This means that, for a triplet satisfying \Eqn{p}, the integers
\begin{subequations}
\begin{eqnarray}
  \trip |a^3-3ab^2| |b^3-3a^2b| c^3 
\end{eqnarray}
form a Pythagorean triplet as well.  The same argument can be applied
to $\tan n\theta$ for any integer $n$, and would mean that the
absolute values of the real and the imaginary parts of $(a+ib)^n$ can
constitute a Pythagorean base if $a$ and $b$ do.
\begin{table}
  \caption{Examples of Pythagorean triplets derived from two such
    triplets.  Note that, when the two triplets are identical, only
    one new triplet results, as should be obvious from
    \Eqn{1and2}.}\label{doubles}
  \begin{center}
    $$
    \begin{array}{c|cccc}
      & \trip 12 5 13 & \trip 8 15 17 & \trip 20 21 29 \\
      \hline
      \multirow2*{$\trip 4 3 5 $} & \trip 16 63 65 &
                 \trip 36 77 85 & \trip 144 17 145 \\
                 & \trip 56 33 65  & \trip 84 13 85 & \trip 24 143 145
                 \\ \hline  
      \multirow2*{$\trip 12 5 13 $} & \multirow2*{$\trip 120 119 169 $} &
      \trip 140 171 221 & \trip 152 345 377 \\ 
                 && \trip 220 21 221 & \trip 352 135 377 \\
    \end{array}
    $$
  \end{center}
\end{table}
For $n=4$ and $=5$, this recipe gives the following Pythagorean
triplets:
\begin{eqnarray}
  & \trip {|4 a^3b - 4 ab^3|} {|b^4 - 6a^2b^2 + a^4|} c^4 ; \\
  & \trip {|a^5 - 10a^3b^2 + 5ab^4|} {|b^5 - 10a^2b^3 + 5a^4b|} 
  c^5 .
\end{eqnarray}
\end{subequations}
Similar expressions can easily be written down with higher powers of
$c$ as the peak.

From this discussion, it is clear that if the peak of a Pythagorean
triplet is expressible in the form $c^n$ for any integer $n$, then the
corresponding triplet can be derived from a triplet whose peak is $c$.
Thus, for example, starting from the triplet $\trip 4 3 5 $. we can
generate not only the triplet $\trip 24 7 25 $ that appears in
\tabl{list}, but also higher triplets like $\trip 44 117 125 $ or
$\trip 336 527 625 $, involving higher powers of 5 as the peak.  These
derivable ones are not fundamental or basic triplets then.

We can identify more which are not basic.  If we have two
Pythagorean triplets $\trip a_1 b_1 c_1 $ and $\trip a_2 b_2 c_2 $
corresponding to the angles $\theta_1$ and $\theta_2$, then the
tangent of the angles $\theta_1+\theta_2$ and $\theta_1-\theta_2$ will
also be rational.  They will generate the new triplets:
\begin{subequations}
  \label{1and2}
\begin{eqnarray}
  & \trip {a_1b_2 + a_2b_1} {|a_1a_2 - b_1b_2|} c_1c_2 \,,
  \label{1+2} \\ 
  & \trip {|a_1b_2 - a_2b_1|} {a_1a_2 + b_1b_2} c_1c_2 \,.
  \label{1-2}
\end{eqnarray}
\end{subequations}
Examples of Pythagorean triplets formed through this recipe have been
shown in \tabl{doubles}.  Notice that when one starts with two copies
of the same triplet, \Eqn{1-2} does not give a valid solution, so
there is only one solution, same as that inferred from
\Eqn{tan2theta}.

This shows two things.  First, any Pythagorean triplet whose peak is a
non-prime number can be generated from other triplets whose peaks are
factors of its own peak.  Second, if the peak is a product of two
prime numbers, it can have two different bases.  Indeed, for peaks of
the form $c_1c_2$, one can find other pairs of numbers whose squared
sum would equal to its square, viz., the pair $(c_1a_2, c_1b_2)$ or
the pair $(a_1c_2, b_1c_2)$.  But these choices do not satisfy
\Eqn{gcd1}, and so we do not consider them at all.

Similar exercises can be performed with three triplets.  If one takes
a product of the type $(a_1\pm ib_1)(a_2 \pm ib_2) \cdots (a_n \pm
ib_n)$ for any number of factors $n$, its real and imaginary parts
will constitute a Pythagorean base if each $a_i$ and $b_i$ do.  For
example, starting from any three triplets, one can construct a triplet
\begin{eqnarray}
  \trip {|a_1a_2a_3 - a_1b_2b_3 - b_1a_2b_3 - b_1b_2a_3|}
  {|b_1a_2a_3 + a_1b_2a_3 + a_1a_2b_3 - b_1b_2b_3|} c_1c_2c_3 \,.
\end{eqnarray}
There will be three more variations, obtained by changing the sign of
one of the $b$'s in this expression.  Thus, for example, if the peak
is $5\times 13\times 17 = 1105$, the solutions for the triplets are
\begin{eqnarray}
  \trip 576 943 1105 , \quad  \trip 1104 47 1105 , \quad
  \trip 264 1073 1105 , \quad \trip 744 817 1105 .
\end{eqnarray}
If there are $n+1$ prime factors in the peak, there can be $2^n$
independent triplets.

The effects can be compounded.  If we consider two angles $\theta_1$
and $\theta_2$ whose tangents are rational numbers, we can construct
the expression of $\tan (n\theta_1+m\theta_2)$, which will contain the
bases of Pythagorean triplets whose peak is $c_1^nc_2^m$.  It is thus
clear that the fundamental triplets have prime numbers as their
Pythagorean peaks.  All others can be constructed from them, by taking
the tangent of suitable combinations of angles.

The question that remains is this: can all primes be Pythagorean
peaks?  Of course 2 and 3 cannot, because they are too small to
accommodate the sum of two squares.  To make a systematic study of the
higher primes, we need to study some properties of the triplets.

\begin{theorem}\label{th:c+-a}

  For any Pythagorian triplet, both $c-a$ and $c+a$ must be perfect
  squares.

  \proof This proof follows simply from \theo{th:allsol}, where we
  showed that there must be two integers $p$ and $q$ such that
  $c=p^2+q^2$ and $a=2pq$.

\end{theorem}

\begin{theorem}
  For any Pythagorean triplet, either $c-b$ or $c+b$ must be divisible
  by 8, but not both.

  \proof We already said the in \Eqn{pq}, $p$ must be even.  Let us
  then write $a=4n_a$, $b=2n_b+1$, $c=2n_c+1$.  Then, $c^2-b^2=a^2$
  reduces to the equation
  \begin{eqnarray}
    (n_c-n_b) (n_c+n_b+1) = 4n_a^2 \,.
  \end{eqnarray}
  The difference between the two factors on the left side is $2n_b+1$,
  which is odd.  This means that one of the factors must be odd and
  the other even.  The even factor must be a multiple of 4 in view of
  the equation.  So there are two options:
  \begin{subequations}
  \begin{eqnarray}
    n_c - n_b = 4k  &\Longrightarrow & c-b = 8k; \\
    n_c+n_b+1 = 4k  &\Longrightarrow & c + b = 8k.
  \end{eqnarray}
  \end{subequations}

\end{theorem}

\begin{theorem}\label{th:cmod4}
  Any Pythagorean peak is of the form $c=4n_c+1$, i.e., $c=1 \mod 4$.

  \proof Since $p=2n_p$ and $q=2n_q+1$ in \Eqn{pq},
  $c=4n_p^2+4n_q^2+4n_q+1 = 1 \mod 4$.

\end{theorem}

\begin{theorem}\label{th:not3}

  A Pythagorean peak cannot be divisible by 3, and one member of the
  base (but not both) must be divisible by 3.

  \proof Any natural number can be written as $n=3k+r_3$, where $r_3$
  is 0 or 1 or 2.  Squaring, one obtains
  \begin{eqnarray}
    n^2 = \options(0 \\ 1) \mod 3 \,,
    \label{sqmod3}
  \end{eqnarray}
  where the numbers in the curly brackets are the possible answers.
  If a Pythagorean peak $c$ is a multiple of 3, so is $c^2$.
  Therefore, $a^2+b^2$ must also be a multiple of 3.  From
  \Eqn{sqmod3}, we see that this is possible only if both $a^2$ and
  $b^2$ are divisible by 3, but in that case, \Eqn{gcd1} is violated,
  so we do not obtain a Pythagorean triplet in the sense we have
  defined them.

  So now that we know that $c$ is not divisible by 3, we have $c^2 = 1
  \mod 3$.  The only way $a^2+b^2$ can also be a number $1\mod3$,
  subject to the result of \Eqn{sqmod3}, is if one of them is
  divisible by 3 and the other leaves a remainder 1 when divided by 3.
  This completes the proof.
  
\end{theorem}

So now we see that there are two kinds of Pythagorean triplets: one
type for which $a$ is a multiple of 3 and $b$ is not, and the other
kind for which these roles are reversed.  Consider the first kind.  We
already showed that $a$ is a multiple of 4.  If it is a also a
multiple of 3, it is a multiple of 12.  Moreover, $b$ cannot be a
multiple of 3, and has to be an odd number.  Therefore, it must be of
the form $6n_b \pm 1$.  For the second kind of triplets, $b$ must be
an odd multiple of 3, i.e., of the form $6n_b+3$, and $a$ should be
divisible by 4 but not by 12, i.e, of the form $4(3n_a \pm1)$.  Here
is a summary:
\begin{subequations}
  \label{2types}
\begin{eqnarray}
  \mbox{\bf Type 1:}   &\qquad&  a = 12n_a,  \quad b = 6n_b \pm1 \,,
  \label{type1} \\
  \mbox{\bf Type 2:}   &\qquad&  a = 4(3n_a \pm 1),  \quad b = 6n_b +3 \,. 
  \label{type2}
\end{eqnarray}
\end{subequations}

\begin{theorem}\label{th:mod24}

  For any Pythagorean triplet,
  \begin{eqnarray}
    c^2 = 1 \mod 24 \,.
  \end{eqnarray}
\end{theorem}
This can be proved by simply squaring the forms of $a$ and $b$ given
in \Eqn{2types} and adding them up.  The following proof might be
simpler in some sense~\cite{24}.
  
  \proof Note that $c^2-1 = (c-1)(c+1)$.  Consider the three
  consecutive integers $c-1$, $c$ and $c+1$.  Since $c$ must be odd by
  \theo{oddc}, both $c-1$ and $c+1$ must be even.  Since they are two
  consecutive even numbers, one of them must be a multiple of 4.
  Thus, one of them is divisible by 2 and the other by 4.  Further,
  among any three consecutive numbers, one must be a multiple of 3.
  Since $c$ cannot be a multiple of 3 through \theo{th:not3}, either
  $c-1$ or $c+1$ is.  Thus, $c^2-1$ is a multiple of
  $2\times4\times3=24$.

\begin{theorem}\label{th:1,5}

  The peak $c$ of a Pythagorean triplet satisfies $c = 1 \mod 12$ or
  $c = 5 \mod 12$.

  \proof We have shown, in \theo{th:cmod4}, that $c = 1 \mod 4$.  In a
  mod-12 arithmetic, it can therefore have the form of either $12k+1$,
  or $12k+5$, or $12k+9$.  But the last one can be written as
  $3(4k+3)$, which is not allowed by \theo{th:not3}.  This leaves us
  with only the other two possibilities, proving the theorem.
  
\end{theorem}

The two kinds of alternatives for the peak given in \theo{th:1,5} are
related to the two types of bases mentioned in \Eqn{2types}.  To see
this, consider first a Type-1 triplet, for which $a$ is a multiple of
12, say $12l$.  If $c=12k+\eta$, then it means that $c \pm a = 12(k
\pm l)+\eta = \eta \mod 12$.  This has to be a perfect square,
according to \theo{th:c+-a}.  However, for any integer $n$,
\begin{eqnarray}
  n^2 = \options(0 \\ 1 \\ 4 \\ 9) \mod 12.
  \label{nsqmod12}
\end{eqnarray}
Since this result applies to $c \pm a$ and $a$ is a multiple of 12, we
must have
\begin{eqnarray}
  c = \options(0 \\ 1 \\ 4 \\ 9) \mod 12. 
\end{eqnarray}
Combining this with the result of \theo{th:1,5}, we find 
\theo{th:a/12} given below.

\begin{theorem}\label{th:a/12}

  For type-1 triplets ($a$ is divisible by 3), $c=1 \mod 12$.

  \proof The proof is given above.

\end{theorem}

Consider now the other case, i.e., we have a Type-2 triplet for which
$a$ is of the form $12l+\eta'$, with $\eta'=4$ or $\eta'=8$.  Then,
\begin{eqnarray}
  c \pm a &=& 12(k-l) + \eta \pm \eta' \,.
\end{eqnarray}
Here is a little list of the possible values of both $\eta + \eta'$
and $\eta-\eta'$ in $\mod 12$ arithmetic, written in that order:
\begin{eqnarray}
  \begin{array}{c|cccc}
    \multirow2*{$\eta'$} & \multicolumn2c \eta \\
    \cline{2-3}
    & 1 & 5 \\
    \hline
    4 & 5, 9 & 9, 1 \\
    8 & 9, 5 & 1, 9 \\ 
  \end{array}
\end{eqnarray}
Comparison with \Eqn{nsqmod12} shows that in this case, only $\eta=5$
is allowed.  This result is summarized in this theorem:

\begin{theorem}

  For type-2 triplets ($b$ is a multiple of 3), $c= 5\mod 12$.

  \proof Already given.
  
\end{theorem}

\begin{table}
  \def\nope#1 {\framebox{#1}}
  \def\ngcd#1 {\begin{picture}(15,15)\put(2,2){#1}
      \put(7,5){\circle{16}} \end{picture}}
  \def\notp#1 {\ensuremath{\bullet}#1\ensuremath{\bullet}}

  \caption{Pythagorean peaks less than 100.}\label{<100} \bigskip
  
    Legend: 
      \nope {$n$} \ means $n$ cannot be a peak, 
      \notp {$n$} \ means $n$ is a non-prime which can be a peak,
      \ngcd {$n$} \ means $n$ can be a peak only if the gcd
      condition is sacrificed.

  $$ 
    \begin{array}{c|cccc}
      k & 12k+1 & 12k+5 & 12k+7 & 12k+11  \\
      \hline
      0 & & 5 & \nope 7 & \nope 11 \\
      1 & 13 & 17 & \nope 19 & \nope 23 \\
      2 & \notp 25 & 29 & \nope 31 & \ngcd 35 \\
      3 & 37 & 41 & \nope 43 & \nope 47 \\
      4 & \nope 49 & 53 & \ngcd 55 & \nope 59 \\
      5 & 61 & \notp 65 & \nope 67 & \nope 71 \\
      6 & 73 & \nope 77 & \nope 79 & \nope 83 \\
      7 & \notp 85 & 89 & \ngcd 91 & \ngcd 95 \\
      8 & 97 &    &    &   \\
    \end{array}
    $$

\end{table}

\begin{theorem}\label{th:not711}

    A Pythagorean peak cannot be divisible by 7 or 11.

    \proof The proof is very similar to that of \theo{th:not3}.  The
    point is that, for every integer $n$,
    \begin{eqnarray}
      n^2 =  \options(0 \\ 1 \\ 2 \\ 4) \mod 7.
      \label{nsqmod7}
    \end{eqnarray}
    If $c=0 \mod7$, we cannot find any nonzero values of $a^2\mod7$
    and $b^2\mod7$ which add up to 7.  That proves the result
    regarding the divisibility by 7.

    Similarly,
    \begin{eqnarray}
      n^2 =  \options(0 \\ 1 \\ 3 \\ 4\\ 5\\ 9) \mod 11.
      \label{nsqmod11}
    \end{eqnarray}
    No pair of two numbers from this list adds up to 11.  Hence the
    result.

\end{theorem}

Results like those in \theo{th:not3} and \theo{th:not711} can be proved
for a great number of primes.  Here is a partial list of such primes:
\begin{eqnarray}
  3,\ 7,\ 11,\ 19,\ 23,\ 31,\ 43,\ 47,\ 59,\ \ldots \,.
\end{eqnarray}
For 3, 7, and 11, we have explicitly shown the proof above.  We can
continue the same kind of proof for the other primes shown here.
Apart from the number 3, all numbers in this list are of the form
$12k+7$ or $12k+11$ for some $k$.  Hence, instead of dealing with each
of them separately, we want to prove the following theorem.

\begin{theorem}

  Any number of the form $12k+7$ or $12k+11$ cannot be a prime factor
  of a Pythagorean peak.

  This says much more than what we have said before.  Earlier, we said
  that Pythagorean peaks cannot be divisible by 7 and 11.  Now we are
  extending the result to every number which leaves a remainder of 7
  or 11 when divided by 12.

  \proof Let
  \begin{eqnarray}
    N = 12k + \options(7 \\ 11).
  \end{eqnarray}
  For either options, it is true that $N = 4l+3$ for some integer
  $l$.  According to \theo{th:cmod4}, we cannot have a Pythagorean
  peak equal to $N$.  The question is, can we have a peak $Nx$ for
  some $x$?

  \theo{th:allsol} tells us that for any even integer $p$ and any odd
  integer $q$, we must have
  \begin{eqnarray}
    p^2 + q^2 \neq N.
    \label{neqN}
  \end{eqnarray}
  If we consider all squares of integers and find their values in
  modulus $N$ arithmetic, we will obtain a list of numbers $r_1 \bsup
  N$, $r_2 \bsup N$ etc, like we got for 7 in \Eqn{nsqmod7} or for 11
  in \Eqn{nsqmod11}.  In other words, we are saying that the square of
  any integer $n$ can be written in the form
  \begin{eqnarray}
    n^2 = N \alpha_n + r_n \bsup N 
    \label{rn}
  \end{eqnarray}
  for some $\alpha_n$, with $0 \leq r_n < N$.  Then \Eqn{neqN} implies
  that, no matter what $p$ and $q$ we choose, we will find that
  \begin{eqnarray}
    r_p \bsup N + r_q \bsup N \neq N \qquad \forall p,q.
    \label{rp+rq}
  \end{eqnarray}
  Now consider a similar list, made from divisions by $Nx$.  We can
  write
  \begin{eqnarray}
    \alpha_n = xy_n + \beta_n ,
  \end{eqnarray}
  which would mean that
  \begin{eqnarray}
    n^2 = N(xy_n + \beta_n) + r_n \bsup N.
  \end{eqnarray}
  If we take out the multiple of $Nx$ from here, we obtain
  \begin{eqnarray}
    r_n \bsup{Nx} =  r_n \bsup N + N \beta_n \,.
  \end{eqnarray}
  If we add two such remainders for $p^2$ and $q^2$, we will get
  \begin{eqnarray}
    r_p \bsup{Nx} +  r_q \bsup{Nx} =  r_p \bsup N + r_q \bsup N + N
    (\beta_p + \beta_q) \,. 
    \label{rNx}
  \end{eqnarray}
  If this number has to equal $Nx$, then it will have to be divisible
  by $N$.  That can happen only if $r_p \bsup N + r_q \bsup N$ is
  divisible by $N$.  $r_p \bsup N + r_q \bsup N=0$ is not allowed
  because that would mean that $r_p \bsup N$ and $r_q \bsup N$ are
  both zero, implying that both $p^2$ and $q^2$ are divisible by $N$,
  which would violate the gcd condition of \Eqn{gcd1}.  \Eqn{rp+rq}
  eliminates the possibility that the sum is equal to $N$.  Any higher
  multiplet is not possible either, since both $r_p \bsup N$ and $r_q
  \bsup N$ are less than $N$.  Thus, the right side of \Eqn{rNx}
  cannot be divisible by $N$, and therefore it cannot be equal to
  $Nx$.  This concludes the proof.

\end{theorem}

In conclusion, then, we find that the only independent Pythagorean
triplets are those whose peak is a prime number, and only primes of
the form $12k+1$ and $12k+5$ are allowed.  All other triplets can be
constructed from these ones by the composition method involving the
trigonometric ratios as described earlier.

As a corollary, we have also shown that if a Pythagorean peak has
$n+1$ prime factors, there can be $2^n$ independent triplets
corresponding to such a peak.  We have also shown that there are two
different types of Pythagorean triplets, with differing
characteristics.

One important question to ask is this: do all primes of the form
$12k+1$ and $12k+5$ qualify as Pythagorean peaks?  I do not know the
answer.  I have not been able to find any that is not.  The search
results in \tabl{<100} do show that all primes of the given forms are
indeed Pythagorean peaks.  But I have not been able to prove that any
such number can be the peak of a Pythagorean triplet.

\paragraph*{Acknowledgements and apologies:} A lot of my motivation
of looking into the Pythagorean triplets came from reading a paper by
my friend Ranjan Ghosh on Fermat's last theorem \cite{ranjan}.  None
of the results presented here directly appears in his paper, except a
few trivial ones like \Eqn{treq} and \theo{th:2odd}, so I do not think
anything here will diminish the value of his paper when he decides to
publish it.  I am also indebted to him, as well as Atish Bagchi, for
comments on an earlier version of this write-up.

I am sure that a lot of these results are known.  I have no formal
training in the subject to know how much is known, or, contrarily,
whether there is any new result in this article.  Even if everything
is known, this article may be treated as my exposition of the known
things.

\section*{Appendix 1}
We have not been able to show that any prime of the form stated in
\theo{th:1,5} can be a Pythagorean peak.  However, we can show that
any multiple of 4 can be used as the even number in the base, and any
odd number (greater than 1) as the odd number in the base.  For a
given $a$, we can take $b=\frac14a^2-1$ and $c=\frac14a^2+1$ to
complete the triplet.  Similarly, for a given odd $b$, we can take
$a=\frac12(b^2-1)$ and $c=\frac12(b^2+1)$, and that would constitute a
triplet.  Clearly, this solution does not work for $b=1$ because it
gives $a=0$, but works for all higher odd numbers.

It is easy to see that all solutions mentioned here satisfy the gcd
criterion, \Eqn{gcd1}.  It is also important to realize that these may
not be the only solution with a given $a$ or a given $b$.

\section*{Appendix 2}
In this appendix, we show that the construction indicated in
\Eqn{1and2} can be inverted.  More explicitly, we have the following
result in mind.

\begin{theorem}

  If there are two Pythagorean triplets of the form $\trip A_1 B_1
  c_1c_2 $ and $\trip A_2 B_2 c_1c_2 $ where $c_1$ and $c_2$ are
  primes, then there also exist triplets with $c_1$ as the peak and
  $c_2$ as the peak.

  \proof All we need to do is the provide a recipe for forming
  separate triplets with $c_1$ and $c_2$ as their peaks.  For this,
  let us start with supposing, without any loss of generality, that
  $A_1>A_2$.  That would also imply $B_1<B_2$, since
  \begin{eqnarray}
    A_1^2 + B_1^2 = A_2^2 + B_2^2 = c_1^2 c_2^2 \,.
    \label{c1c2sq}
  \end{eqnarray}
  Now define the following numbers:
  \begin{subequations}
    \label{2gcds}
  \begin{eqnarray}
    \gcd(A_1+A_2, B_2-B_1) &=& 2b_1 \,,
    \label{2gcd1} \\
    \gcd(A_1-A_2, B_2-B_1) &=& 2b_2 \,. 
    \label{2gcd2}
  \end{eqnarray}
  \end{subequations}
  Notice that all combinations of the $A$'s and the $B$'s appearing
  here are even, so there must be a common factor of 2 between the
  pairs, which has been separated out.  Also, note that the
  combination $B_2+B_1$ does not appear in these equations, because it
  is not an independent combination, in view of the fact that
  \begin{eqnarray}
    A_1^2 - A_2^2 = B_2^2 - B_1^2
    \label{ABreln}
  \end{eqnarray}
  which follows from \Eqn{c1c2sq}.
  
  Now, \Eqn{2gcd1} implies that we can write
  \begin{subequations}
    \label{ABeqs}
  \begin{eqnarray}
    A_1+A_2 = 2 a_2b_1, \qquad B_2-B_1 = 2 b_2' b_1
  \end{eqnarray}
  for suitable numbers $a_2$ and $b_2'$.  Similarly, using \Eqn{2gcd2}
  we can write
  \begin{eqnarray}
    A_1-A_2 = 2a_1b_2, \qquad B_2-B_1 = 2b_1'b_2 \,.
  \end{eqnarray}
  \end{subequations}
  Consistency of these two equations demand that
  \begin{eqnarray}
    b_1' = b_1, \qquad b_2' = b_2 \,.
  \end{eqnarray}
  Putting these back into \Eqn{ABeqs} and using \Eqn{ABreln}, we find 
  \begin{eqnarray}
    A_1 &=& a_2b_1 + a_1b_2 \,, \nonumber\\
    A_2 &=& a_2b_1 - a_1b_2 \,, \nonumber\\
    B_1 &=& a_1a_2 - b_1b_2 \,, \nonumber\\
    B_2 &=& a_1a_2 + b_1b_2 \,
  \end{eqnarray}
  Putting these expressions into \Eqn{c1c2sq}, we obtain
  \begin{eqnarray}
    (a_1^2 + b_1^2) (a_2^2 + b_2^2) = c_1^2 c_2^2 \,.
  \end{eqnarray}
  Certainly, this equation is satisfied with either $c_1^2$ or $c_2^2$
  equal to $a_1^2+b_1^2$, and the other one equal to $a_2^2+b_2^2$.
  This completes the proof of the theorem.

\end{theorem}

\end{document}